\documentclass[12pt]{article}

\usepackage{theorem,amssymb,amsmath}

\advance \topmargin by -\headheight
\advance \topmargin by -\headsep
\evensidemargin \oddsidemargin
\author{J.-P. Allouche\thanks{The author was partially supported by the ANR
project ``FAN'' (Fractals et Num\'eration).} \\
CNRS, Institut de Math. de Jussieu \\
Universit\'e P. et M. Curie, Case 247 \\
4 Place Jussieu \\
F-75252 Paris Cedex 05, France \\
{\tt allouche@math.jussieu.fr}
\and
Benoit Cloitre \\
19, rue Louise Michel \\
92300 Levallois-Perret France \\
{\tt benoit7848c@yahoo.fr}
\and
V. Shevelev \\
Department of Mathematics \\
Ben-Gurion University of the Negev \\
Beersheva, Israel \\
{\tt shevelev@bgu.ac.il}
}

\title{Beyond odious and evil}

\date{ }

\def \proof{\bigbreak\noindent{\it Proof.\ \ }}

\def \endpf{{\ \ $\Box$ \medbreak}}

\newtheorem{theorem}{Theorem}
\newtheorem{lemma}{Lemma}
\newtheorem{corollary}{Corollary}
\newtheorem{proposition}{Proposition}

\theorembodyfont{\rm}

\newtheorem{remark}{Remark}
\newtheorem{example}{Example}
\newtheorem{definition}{Definition}

\begin{document}

\maketitle

\hskip 3truecm [...] The lights hanging 

\hskip 3truecm from oak beams above the readers

\hskip 3truecm light {\it and} illuminate every page. Each book dusted

\hskip 3truecm each day. Original jackets, no odious numbers glued to spines,

\hskip 3truecm not one decimal, Dewey or otherwise, in the entire place! [...]

\hskip 3truecm (Thomas Lux, {\it The Ambrosiana Library})

\begin{abstract}
In a recent post on the Seqfan list the third author proposed a conjecture concerning the 
summatory function of odious numbers (i.e., of numbers whose sum of binary digits is odd),
and its analog for evil numbers (i.e., of numbers whose sum of binary digits is even).
We prove these conjectures here. We will also study the sequences of ``generalized'' 
odious and evil numbers, and their iterations, giving in particular a characterization
of the sequences of usual odious and evil numbers in terms of functional equations satisfied
by their compositions.

\medskip

\noindent
{\bf Keywords}: Odious numbers; odd numbers; Thue-Morse sequence; summatory functions; iteration of sequences.

\medskip

\noindent
{\bf MSC Classes}: 11A63, 11B83, 11B85, 11A07, 05A15.

\end{abstract}

\section{Introduction}

The purpose of this paper, whose title is reminiscent of \cite{Nietzsche}, is revisiting the study 
of two families of integers respectively called odious and evil numbers, as well as the study of some
generalizations. A natural integer is called {\it odious} if the sum of its binary digits is odd. 
A natural integer is called {\it evil} if the sum of its binary digits is even. This terminology 
was introduced by the authors of \cite{BCK, BCK2}, see \cite[p.\ 463]{BCK2}; the words ``odious'' 
and ``evil'' were chosen because they begin respectively like ``odd'' and ``even''. 
Let  ${\mathbf a} = (a(n))_{n \geq 0}$ denote the increasing sequence of odious numbers, and  
${\mathbf b} = (b(n))_{n \geq 0}$ denote the increasing sequence of evil numbers. Sequences 
${\mathbf a}$ and ${\mathbf b}$ are respectively A000069 and A001969 in \cite{oeis}, except that 
we let the indexes start from $0$ instead from $1$. Sequences ${\mathbf a}$ and ${\mathbf b}$ begin
$$
{\mathbf a} = 1 \ 2 \ 4 \ 7 \ 8 \ 11 \ 13 \ 14 \ 16 \ 19 \ 21 \ 22 \ 25 \ 26 \ 28 \ 31 \ \ldots
$$
$$
{\mathbf b} = 0 \ 3 \ 5 \ 6 \ 9 \ 10 \ 12 \ 15 \ 17 \ 18 \ 20 \ 23 \ 24 \ 27 \ 29 \ 30 \ \ldots
$$

\begin{remark}
We seem to remember having read somewhere on the web (but where was it?) that mathematicians 
probably do not like numbers, since for them numbers are necessarily either evil or they are odd!
Also note that for some authors the expression ``evil numbers'' has a quite different meaning. 
This is explained, e.g., in \cite{Wolfram}: a number is called evil in that terminology if the 
first $n$ decimal digits of its fractional part sum to $666$ for some integer $n$. The expression 
``evil numbers'' is also used with other meanings, inspiring in particular artists, like 
Fabio Mauri (see \cite{Emmer}).
\end{remark}

The purpose of this paper is twofold. First to prove the above conjecture. Second to describe 
iterations of sequences ${\mathbf a}$ and ${\mathbf b}$ and generalizations. Namely it was noted by
the second author in \cite[A000069]{oeis} that $a(a(n)) = 2 a(n)$.
We will generalize this result by proving Theorem~\ref{iter} below which gives an expression of
$a_{j,d}(a_{i,d}(n))$, where $(a_{j,d}(n))_{n \geq 0}$ denote the increasing sequence of integers 
whose sum of $d$-ary digits is congruent to $j$ modulo $d$.

\section{Iteration of the sequences of generalized odious and evil numbers}

As recalled in the introduction, it was noted by the second author \cite[A000069]{oeis} that the 
increasing sequence of odious numbers ${\mathbf a}$ satisfies $a(a(n)) = 2 a(n)$.
We will generalize this result by proving Theorem~\ref{iter} below.

\bigskip

\noindent
We begin with a definition and a lemma. 

\begin{definition}\label{def}

\ { }

\begin{itemize}
\item Let $d \geq 2$ be an integer. For any integer $x$ we let $\overline{(x)}_d$ denote the 
      residue modulo $d$ of $x$, i.e., the integer belonging to $[0, d-1]$ and congruent to 
      $x$ modulo $d$. Note that we have $x = d \lfloor\frac{x}{d}\rfloor + \overline{(x)}_d$.
\item We let $s_d(n)$ denote the sum of the $d$-ary digits of $n$. We let 
      ${\mathbf t_d} = (t_d(n))_{n \geq 0}$ denote the sequence of integers defined by 
      $t_d(n) \equiv s_d(n) \bmod d$ and $0 \leq t_d(n) \leq d-1$, i.e., 
      $t_d(n) = \overline{(s_d(n))}_d$.
\item For $j \in [0, d-1]$, we let ${\mathbf a_{j,d}} = (a_{j,d}(n))_{n \geq 0}$ denote the
      increasing sequence of integers $k$ such that $s_d(k) \equiv j \bmod d$ (i.e., $t_d(k) = j$).
\end{itemize}
\end{definition}

\begin{lemma}\label{td-morphic}

 \ \ { } \

\begin{itemize}

\item
Sequence ${\mathbf t_d}$ is the fixed point of the morphism defined on $\{0, 1, \ldots, d-1\}$ by 
$0 \to 0 \ \ 1 \ \ \ldots \ \ d-1$, $1 \to 1 \ \ 2 \ \ \ldots \ \ d-1 \ 0$,...,
$d-1 \to d-1 \ \ 0 \ \ 1 \ \ \ldots \ \ d-2$.

\item
If $\alpha$ belongs to $[0, d-1]$, then, for all $n \geq 0$, we have
$$
t_d(dn + \alpha) = \overline{(t_d(n)+\alpha)}_d =
\begin{cases}
t_d(n) + \alpha    \ &\mbox{\rm if $t_d(n) + \alpha \leq d-1$} \\
t_d(n) + \alpha -d \ &\mbox{\rm if $t_d(n) + \alpha \geq d$}. \\
\end{cases}
$$

\item 
If $j$ belongs to $[0, d-1]$, then, for all $n \geq 0$, we have
$$
d - 1 - t_d(dn+d-1-j) =  \begin{cases} j - t_d(n)  &\mbox{\rm if} \ 0 \leq t_d(n) \leq j \\
                             d + j - t_d(n) &\mbox{\rm if} \ j+1 \leq t_d(n) \leq d-1.
                  \end{cases}
$$

\end{itemize}

\end{lemma}

\proof The proof of the first two items is easy and left to the reader. The last item is an
easy consequence of the second item. \endpf

\bigskip

Now we prove a helpful proposition.

\begin{proposition}\label{helpful} 
The sequence $(a_{j,d}(n))_{n \geq 0}$ satisfies
$$
a_{j,d}(n) = dn + \overline{(j-t_d(n))}_d.
$$
This can also be written
$$
\begin{array}{lll}
a_{j,d}(n) &=& dn + \begin{cases} j - t_d(n)  &\mbox{\rm if} \ 0 \leq t_d(n) \leq j \\
                             d + j - t_d(n) &\mbox{\rm if} \ j+1 \leq t_d(n) \leq d-1
                  \end{cases} \\
           & & \\
           &=& dn + d - 1 - t_d(dn + d - i - 1).
\end{array}
$$
\end{proposition}

\proof These equalities are easy consequences of Lemma~\ref{td-morphic}, which implies 
in particular that sequence ${\mathbf t_d}$ consists of consecutive blocks taken from
$(0 \ \ 1 \ \ \ldots \ \ d-1)$, $(1 \ \ 2 \ \ \ldots \ \ d-1 \ \ 0)$,...,
$(d-1 \ \ 0 \ \ 1 \ \ \ldots \ \ d-2)$, where the $r$-th block begins with $t_d(r)$.   \endpf

\bigskip

We are ready to state and prove the result of this section.

\begin{theorem}\label{iter}
For all $n \geq 0$, for all $i, j \in [0, d-1]$, we have
$$
a_{j,d}(a_{i,d}(n)) =  d a_{i,d}(n) + \overline{(j-i)}_d =
\begin{cases}
d a_{i,d}(n) + j - i     \ &\mbox{\rm if $j \geq i$} \\
d a_{i,d}(n) + d + j - i \ &\mbox{\rm if $j < i$}
\end{cases}
$$
\end{theorem}

\proof  Using Proposition~\ref{helpful} we can write 
$$
a_{j,d}(a_{i,d}(n)) = da_{i,d}(n) + \overline{j-t_d(a_{i,d}(n))}_d.
$$
But, from the definition of $a_{i,d}$ we have $t_d(a_{i,d}(n)) = i$. Hence
$$
a_{j,d}(a_{i,d}(n)) = da_{i,d}(n) + \overline{(j-i)}_d =
\begin{cases}
d a_{i,d}(n) + j - i     \ &\mbox{\rm if $j \geq i$} \\
d a_{i,d}(n) + d + j - i \ &\mbox{\rm if $j < i$}. \ \ \Box
\end{cases}
$$

\section{Summatory function of generalized odious and evil numbers}

In order to address Shevelev's conjecture recalled in the introduction, we have to study
the summatory function of odious numbers. This section is devoted to studying the summatory
function of generalized odious and evil numbers. The first step is the following proposition.
(We keep the notation in Definition~\ref{def}.)

\begin{proposition}\label{prop-sum}
Let $a$ and $r$ be integers in $[0, d-1]$. Then
$$
\sum_{\ell = 0}^r \overline{(a - \ell)}_d =
\begin{cases}
a(r+1) - \displaystyle\frac{r(r+1)}{2} &\mbox{\rm if } r \leq a \\
a(r+1-d) + dr - \displaystyle\frac{r(r+1)}{2} &\mbox{\rm if } r > a.
\end{cases}
$$
This can also be written
$$
\sum_{\ell = 0}^r \overline{(a - \ell)}_d =
a(r+1) - \frac{r(r+1)}{2} + d\max\{r-a, 0\}.
$$
\end{proposition}

\proof If $r \leq a$, then
$$
\begin{array}{lll}
\displaystyle\sum_{\ell = 0}^r \overline{(a - \ell)}_d 
&=& \overline{(a)}_d + \overline{(a-1)}_d + \cdots + \overline{(a-r)}_d \\
&=& a + (a-1) + \cdots + (a-r) = a(r+1) - \displaystyle\frac{r(r+1)}{2}\cdot \\
\end{array}
$$

Now, if $r > a$, then
$$
\begin{array}{lll}
\displaystyle\sum_{\ell = 0}^r \overline{(a - \ell)}_d 
&=& \overline{(a)}_d + \overline{(a-1)}_d + \cdots + 
    \overline{(0)}_d + \overline{(-1)}_d + \overline{-(r-a)}_d \\
&=& a + (a-1) + \cdots + 0 + (d-1) + \cdots + (d-r+a) \\
&=& a(r+1-d) + dr - \displaystyle\frac{r(r+1)}{2}\cdot   \ \Box
\end{array}
$$

\bigskip
\bigskip

\noindent
Using Proposition~\ref{helpful} and Proposition~\ref{prop-sum} we obtain
the following theorem.

\begin{theorem}\label{th-sum}
The summatory function of the sequence ${\mathbf a}_{j,d}$ is given by
\begin{multline*}
\sum_{k=0}^N a_{j,d}(k) = \frac{dN(N+1)}{2} + \frac{\lfloor N/d \rfloor d(d-1)}{2} + 
\overline{(j - t_d(\lfloor N/d \rfloor))}_d (\overline{(N)}_d + 1) \\
         - \frac{\overline{(N)}_d(\overline{(N)}_d + 1)}{2} 
+ d\max\{\overline{(N)}_d - \overline{(j - t_d(\lfloor N/d \rfloor))}_d, 0\}
\end{multline*}
\end{theorem}

\proof  We first note that, for any integer $a$, we have
$$
\sum_{k=0}^{ad - 1} \overline{(j-t_d(k))}_d = 
   \sum_{j=0}^{a-1}\sum_{k=jd}^{(j+1)d - 1} \overline{(j-t_d(k))}_d = 
   \sum_{j=0}^{a-1}\frac{d(d-1)}{2} = \frac{ad(d-1)}{2}
$$
since for each $j$, when $k$ runs in $[jd, (j+1)d - 1]$, $t_d(k)$ takes exactly once every value 
in $[0, d-1]$, thus $\overline{(j-t_d(k))}_d$ also takes exactly once every value in $[0,d-1]$.

\medskip

Now, using Proposition~\ref{helpful}, we get
$$
\begin{array}{lll}
\displaystyle\sum_{k=0}^N a_{j,d}(k) 
  &=& \displaystyle\sum_{k=0}^N (dk + \overline{(j-t_d(k))}_d)
          = \displaystyle\frac{dN(N+1)}{2} + \displaystyle\sum_{k=0}^N \overline{(j-t_d(k))}_d \\
  &=& \displaystyle\frac{dN(N+1)}{2} + 
         \displaystyle\sum_{k=0}^{d\lfloor N/d \rfloor - 1} \overline{(j-t_d(k))}_d +
         \displaystyle\sum_{k=d\lfloor N/d \rfloor}^N \overline{(j-t_d(k))}_d \\
  &=& \displaystyle\frac{dN(N+1)}{2} +
      \displaystyle\frac{\lfloor N/d \rfloor d (d-1)}{2} +
      \displaystyle\sum_{k=d\lfloor N/d \rfloor}^N \overline{(j-t_d(k))}_d.
\end{array}
$$
But
\begin{multline*}
\sum_{k=d\lfloor N/d \rfloor}^N \overline{(j-t_d(k))}_d =
\sum_{\ell=0}^{N-d\lfloor N/d \rfloor} \overline{(j-t_d(\ell + d\lfloor N/d \rfloor))}_d = \\
\sum_{\ell=0}^{N-d\lfloor N/d \rfloor} \overline{(j-\overline{(t_d(\ell) + t_d(\lfloor N/d \rfloor)}_d)}_d
= \sum_{\ell=0}^{N-d\lfloor N/d \rfloor} \overline{(j-t_d(\lfloor N/d \rfloor) - \ell)}_d = \\
\sum_{\ell=0}^{N-d\lfloor N/d \rfloor} \overline{(\overline{(j-t_d(\lfloor N/d \rfloor))}_d - \ell)}_d. 
\end{multline*}
Now, using Proposition~\ref{prop-sum} with $a$ replaced by $\overline{(j-t_d(\lfloor N/d \rfloor))}_d$
and $r$ replaced by $N-d\lfloor N/d \rfloor$ yields the result. \endpf

\bigskip

Now we give a corollary of Theorem~\ref{th-sum} above, recalling in passing some notation 
of the introduction.

\begin{corollary}\label{easy}
Let ${\mathbf t} = (t(n))_{n \geq 0}$ be the Thue-Morse sequence defined by $t(0) = 0$, and
for all $n \geq 0$, $t(2n) = t(n)$, $t(2n+1) = 1 - t(n)$.
Let ${\mathbf a} = (a(n))_{n \geq 0}$ be the sequence of odious numbers.
Let ${\mathbf b} = (b(n))_{n \geq 0}$ be the sequence of evil numbers.
Let $S(n) = a(0) + a(1) + \cdots + a(n)$ be the summatory function of sequence ${\mathbf a}$.
Let $R(n) = b(0) + b(1) + \cdots + b(n)$ be the summatory function of sequence ${\mathbf b}$.
Then

\begin{itemize}

\item $a(n) = 2n + 1 - t(n)$

\item $b(n) = 2n + t(n)$

\item $S(n) = n^2 + \frac{3n}{2} + \frac{1}{2} + \frac{1+(-1)^n}{4} (1 - 2t(n))
= \left\{
\begin{array}{ll}
n^2 + \frac{3n}{2} + \frac{1}{2} \ &\mbox{\rm if $n$ is odd} \\
n^2 + \frac{3n}{2} + 1 - t(n) \ &\mbox{\rm if $n$ is even.}
\end{array}
\right.
$

\item $R(n) = n^2 + \frac{3n}{2} + \frac{1}{2} + \frac{1+(-1)^n}{4}(2t(n)-1)
= \left\{
\begin{array}{ll}
n^2 + \frac{3n}{2} + \frac{1}{2} \ &\mbox{\rm if $n$ is odd} \\
n^2 + \frac{3n}{2} + t(n) \ &\mbox{\rm if $n$ is even}.
\end{array}
\right.
$
\end{itemize}

\end{corollary}

\proof Put $d=2$ and $j=0,1$ in Theorem~\ref{th-sum}. (Also see, e.g., \cite[Section~8]{guy},
\cite[A000069]{oeis} and \cite[A173209]{oeis}.)

\bigskip

Note that it is also possible to give summatory functions of polynomial expressions for sequences 
like ${\mathbf a}$ and ${\mathbf b}$. For example, we can prove the following result.

\begin{corollary}
We have the relations
$$
\begin{array}{lll}

\displaystyle\sum_{k=0}^{b(n)} a(k) &=& b(n)^2+ b(n) + n + 1 \\

\displaystyle\sum_{k=0}^{a(n)} b(k) &=& a(n)^2 + a(n) + n + 1 \\

\displaystyle\sum_{k=0}^{a(n)} a(k) &=& a(n)^2 + 2a(n) - n \\

\displaystyle\sum_{k=0}^{b(n)} b(k) &=& b(n)^2 + 2b(n) - n^2. \\

\end{array}
$$
We have the relations
$$
\sum_{a(k) \leq n} a(k) = \left\{
\begin{array}{ll}
\frac{n^2}{4} - \frac{n}{4} + nt(n) \ &\mbox{\rm if $n \equiv 0 \bmod 4$} \\
\frac{n^2}{4} + \frac{n}{4} - \frac{1}{2} + t(n) \ &\mbox{\rm if $n \equiv 1 \bmod 4$} \\
\frac{n^2}{4} - \frac{n}{4} - \frac{1}{2} + (n+1)t(n) \ &\mbox{\rm if $n \equiv 2 \bmod 4$} \\
\frac{n^2}{4} + \frac{n}{4} \ &\mbox{\rm if $n \equiv 3 \bmod 4$} \\
\end{array}
\right.
$$
$$
\sum_{b(k) \leq n} b(k) = \left\{
\begin{array}{ll}
\frac{n^2}{4} + \frac{3n}{4} - nt(n) \ &\mbox{\rm if $n \equiv 0 \bmod 4$} \\
\frac{n^2}{4} + \frac{n}{4} + \frac{1}{2} - t(n) \ &\mbox{\rm if $n \equiv 1 \bmod 4$} \\
\frac{n^2}{4} + \frac{3n}{4} + \frac{1}{2} - (n+1)t(n) \ &\mbox{\rm if $n \equiv 2 \bmod 4$} \\
\frac{n^2}{4} + \frac{n}{4} \ &\mbox{\rm if $n \equiv 3 \bmod 4$}. \\
\end{array}
\right.
$$
\end{corollary}

\proof Left to the reader. Hint for the last two formulas: note that
$$
\begin{array}{lll}
\{k; \ a(k) \leq n\} &=& \{k; \ 2k+1-t(k) \leq n\} \\
&=& \{k; \ t(k)=1, \ k \leq \lfloor \frac{n}{2} \rfloor \} 
\cup \{k; \ t(k)=0, \ k \leq \lfloor \frac{n-1}{2}\rfloor\} \\
&=& \{0, 1, 2, \ldots, \lfloor \frac{n-1}{2} \rfloor \} \cup \Theta_n
\end{array}
$$
where
$$
\Theta_n =
\left\{
\begin{array}{ll}
          \ \ \ \emptyset \ &\mbox{\rm if $n$ is odd, or if $n$ is even and $t(n)=0$} \\
          \{\lfloor \frac{n}{2}\rfloor\} 
          \ &\mbox{\rm if $n$ is even and $t(n)=1$.}
\end{array}
\right.
$$
Thus
$$
\sum_{a(k) \leq n} a(k) = \sum_{k \leq \lfloor \frac{n-1}{2} \rfloor} a(k) 
+ \frac{1+(-1)^n}{2} t(n) a(\lfloor n/2 \rfloor). \ \Box
$$

\section{Proof of a conjecture of Shevelev}

Before stating Shevelev's conjecture, we need a definition and a lemma.

\begin{definition}\label{mod4}
If $x$ and $y$ are two integers, we write $x <_4 y$ (resp.\ $x \leq_4 y$) if the residues
modulo $4$ of $x$ and $y$, denoted by $\overline{x}$ and $\overline{y}$, belonging to
$\{0, 1, 2, 3\}$ and considered as natural integers, satisfy $\overline{x} < \overline{y}$
(resp.\ $\overline{x} \leq \overline{y}$).
\end{definition}

\begin{example}
For example $17 <_4 6$ because $17 \equiv 1 \bmod 4$, $6 \equiv 2 \bmod 4$ and $1 < 2$.
\end{example}

\begin{lemma}\label{ineq}
Let ${\mathbf a}$ and ${\mathbf t}$ be as above the increasing sequence of odious numbers
and the Thue-Morse sequence. Then, if $n$ and $m$ are both odd or both even, then
$a(n) <_4 a(m)$ (resp.\ $a(n) \leq_4 a(m)$) if and only if $t(m) < t(n)$ (resp.\ $t(m) \leq t(n)$).
\end{lemma}

\proof Since $a(n) = 2n + 1 - t(n)$, we see that $\overline{a(n)} = 1 - t(n)$ if $n$ is even,
and that $\overline{a(n)} = 3 - t(n)$ if $n$ is odd. The statement in the lemma follows. \endpf

\bigskip

\begin{theorem}[Shevelev's conjecture]\label{conj}
Let $S(n)$ be the summatory function of odious numbers, i.e., $S(n) = a(0) + a(1) + \cdots + a(n)$.
We have, for $n \geq 2$,
\begin{itemize}

\item if $a(n-1) <_4 a(n+1)$ and $a(n) \leq_4 a(n+2)$,
      then $S(n) = \displaystyle\frac{a(n)a(n+1)}{4}$

\item if $a(n-1) >_4 a(n+1)$ and $a(n) \geq_4 a(n+2)$,
      then $S(n) = \displaystyle\frac{a(n)a(n+1)}{4} + \frac{1}{2}$

\item if $a(n-1) <_4 a(n+1)$ and $a(n) >_4 a(n+2)$,
      then $S(n) = \displaystyle\frac{a(n)(a(n+1)-1)}{4}$

\item if $a(n-1) >_4 a(n+1)$ and $a(n) <_4 a(n+2)$,
      then $S(n) = \displaystyle\frac{(a(n)+1)a(n+1)}{4}\cdot$

\end{itemize}
\end{theorem}

\proof 

\begin{itemize}

\item[(i)] Proof of the first assertion. From Lemma~\ref{ineq} the conditions $a(n-1) <_4 a(n+1)$
and $a(n) \leq_4 a(n+2)$ are equivalent to $t(n+1) < t(n-1)$ and $t(n+2) \leq t(n)$.
Thus $t(n+1) = 0$, $t(n-1) = 1$, and either $t(n) = 1$ or $t(n) = t(n+2) = 0$.
But $t(n) = t(n+2) = 0$ is impossible because this would imply $((t(n), t(n+1), t(n+2)) = (0, 0, 0)$
and the Thue-Morse sequence does not contain any cube. Thus $t(n) = 1$, $t(n-1) = 1$, $t(n+1) = 0$.
Furthermore $t(n) = t(n-1) (= 1)$ implies that $n$ must be even (if $n = 2k+1$, $t(n) = 1-t(k)$ and
$t(n-1) = t(k)$). So, using Corollary~\ref{easy}, 
$S(n) = n^2 + \frac{3n}{2} + 1 - t(n) = n^2 + \frac{3n}{2}$. On the other hand
$a(n)a(n+1) = (2n+1-t(n))(2n+3-t(n+1) = 2n(2n+3) = 4n^2 + 6n = 4(n^2 + \frac{3n}{2})$.

\item[(ii)] Proof of the second assertion. From Lemma~\ref{ineq} the conditions $a(n-1) >_4 a(n+1)$
and $a(n) \geq_4 a(n+2)$ are equivalent to $t(n+1) > t(n-1)$ and $t(n+2) \geq t(n)$. Hence $t(n-1) = 0$,
$t(n+1) = 1$, and either $t(n) = t(n+2) = 1$ or $t(n) = 0$. But we cannot have $t(n) = t(n+2) = 1$, because
this would give $(t(n), t(n+1), t(n+2)) = (1, 1, 1)$ and this would give a cube in the Thue-Morse sequence.
Thus $t(n-1) = 0$, $t(n) = 0$, $t(n+1) = 1$. As previously $t(n-1) = t(n) (= 0)$ implies that $n$ must be
even. So, using Corollary~\ref{easy}, $S(n) = n^2 + \frac{3n}{2} + 1 - t(n) = n^2 + \frac{3n}{2} + 1$. 
On the other hand
$a(n)a(n+1) + 2 = (2n+1-t(n))(2n+3-t(n+1) + 2 = (2n+1)(2n+2) + 2 = 4n^2 + 6n + 4 = 4(n^2 + \frac{3n}{2} + 1)$.

\item[(iii)] Proof of the third assertion. From Lemma~\ref{ineq} the conditions $a(n-1) <_4 a(n+1)$
and $a(n) >_4 a(n+2)$ are equivalent to $t(n+1) < t(n-1)$ and $t(n+2) > t(n)$. Thus $t(n-1) = 1$,
$t(n) = 0$, $t(n+1) = 0$, $t(n+2) = 1$. Since $t(n) = t(n+1) (= 0)$, $n$ must be odd. So, using 
Corollary~\ref{easy}, $S(n) = n^2 + \frac{3n}{2} + \frac{1}{2}$. On the other hand $a(n)(a(n+1) - 1) = 
(2n + 1 - t(n))(2n + 2 - t(n+1)) = (2n+1)(2n+2) = 4(n^2 + \frac{3n}{2} + \frac{1}{2})$.

\item[(iv)] Proof of the fourth assertion. From Lemma~\ref{ineq} the conditions $a(n-1) >_4 a(n+1)$
and $a(n) <_4 a(n+2)$ are equivalent to $t(n+1) > t(n-1)$ and $t(n+2) < t(n)$. Hence $t(n-1) = 0$,
$t(n) = 1$, $t(n+1) = 1$, $t(n+2) = 0$. Since $t(n) = t(n+1) (= 1)$, $n$ must be odd. So, using 
Corollary~\ref{easy}, $S(n) = n^2 + \frac{3n}{2} + \frac{1}{2}$. On the other hand 
$(a(n)+1)a(n+1) = (2n+2-t(n))(2n+3-t(n+1)) = (2n+1)(2n+2) = 4(n^2 + \frac{3n}{2} + \frac{1}{2})$. \endpf

\end{itemize}

\section{Functional equations for sequences ${\mathbf a}$ and ${\mathbf b}$}

Several studies about iterating increasing sequences of integers can be found in the literature
(see, e.g., \cite{ARS, CSV, LY, Sarkaria} and references therein, in particular parts of 
Reference \cite{Targonski} that we discovered thanks to \cite{Sarkaria}).

With the previous notation, the increasing sequences of odious and of evil numbers
satisfy $a(n) = a_{1,2}(n)$ and $b(n) = a_{0,2}(n)$. We thus have the following
relations.

\begin{corollary}\label{iterations}

\ { } 

\begin{itemize}

\item{{\rm (i)}} $a(a(n)) = 2a(n)$

\item{{\rm (ii)}} $b(b(n)) = 2b(n)$

\item{{\rm (iii)}} $a(b(n)) = 2b(n) + 1$

\item{{\rm (iv)}} $b(a(n)) = 2a(n) + 1$

\item{{\rm (v)}} $a(a(n)) = b(a(n)) - 1$

\item {{\rm (vi)}} $b(b(n)) = a(b(n)) - 1$

\item{{\rm (vii)}} $a(n) - b(n) = 1 - 2 t(n)$ 
                   (in particular $a(n) - b(n)$ takes only the values $\pm 1$)

\item{{\rm (viii)}} $a(b(n)) - b(a(n)) = 4t(n) - 2$ 
                    (in particular $a(b(n)) - b(a(n))$ takes only the values $\pm 2$).

\end{itemize}

\end{corollary}

\proof The first four relations are Theorem~\ref{iter} for the case $d=2$. 
Relations~(v) and (vi) are easy consequences of Relations~(i) to (iv).
The last two relations are consequences of the expressions of $a(n)$ and $b(n)$
given in Corollary~\ref{easy} and of the properties $t(2n) = t(n)$ and 
$t(2n+1) = 1 - t(n)$. \endpf

\bigskip

One might ask which set of relations among relations (i) to (vi) suffice to
characterize sequences ${\mathbf a}$ and ${\mathbf b}$. The next three theorems
yield three answers to the question.

\begin{theorem}~\label{equat}
Suppose that the two sets $X$ and $Y$ form a partition of the integers. 
Let ${\mathbf x} = (x_n)_{n \geq 0}$ be the increasing sequence of the elements of $X$,
and let ${\mathbf y} = (y_n)_{n \geq 0}$ be the increasing sequence of the elements of $Y$.
Suppose that ${\mathbf x}$ and ${\mathbf y}$ satisfy the following relations

\begin{itemize}

\item $x(x(n)) = 2x(n)$ for all $n \geq 0$

\item $y(y(n)) = 2y(n)$ for all $n \geq 0$

\item $|x(n) - y(n)| = 1$ for all $n \geq 0$
\end{itemize}

Then, either ${\mathbf x} = {\mathbf a}$ and ${\mathbf y} = {\mathbf b}$, or
${\mathbf x} = {\mathbf b}$ and ${\mathbf y} = {\mathbf a}$. In particular
the sequence $(x(n) - y(n))_{n \geq 0}$ must be equal to $(1 - 2t(n))_{n \geq 0}$
or to $(2t(n) - 1)_{n \geq 0}$.
\end{theorem}

\proof We must have that $\{0, 1\} = \{x(0), y(0)\}$. Without loss of generality 
we may suppose that $x(0) = 1$ thus $y(0) = 0$. We thus want to prove that 
${\mathbf x} = {\mathbf a}$ and ${\mathbf y} = {\mathbf b}$.
We will prove by induction on $n$ that $\{2n, 2n+1\} = \{x(n), y(n)\}$. The property is true for $n=0$;
suppose it is true for $n$ and let us look at $\{2n+2, 2n+3\}$. Either there exists $k$ such that
$2n+2 = x(k)$ or there exists $k$ such that $2n+2 = y(k)$ ($X$ and $Y$ form a partition of the integers).

\medskip

If $2n+2 = x(k)$ we have necessarily $2n+3 = y(k)$ (since $|x(k) - y(k)| = 1$). Furthermore $k \geq n+1$
(since ${\mathbf x}$ and ${\mathbf y}$ are increasing). If we had $k \geq n+2$ the values $2n+2$, $2n+3$
would not be in the range of ${\mathbf x}$ nor in the range of ${\mathbf y}$, hence $k = n+1$.

\medskip

If $2n+2 = y(k)$, the same reasoning shows that $2n+3 = x(k)$, and $k = n+1$.

\medskip

\noindent
We thus have $\{2n+1, 2n+3\} = \{x(n+1), y(n+1)\}$ and the induction is proven. Now, define the 
sequence $(\alpha(n))_{n \geq 0}$ by $x(n) = 2n + 1 - \alpha(n)$. This implies of course $\alpha(0)=0$ 
and $y(n) = 2n + \alpha(n)$. We then note that, for any integer $m$, we have, by applying the formula
$x(n) = 2n + 1 - \alpha(n)$ with $n=x(m)$, on one hand $x(x(m)) = 2x(m) + 1 - \alpha(x(m))$, 
and on the other hand $x(x(m)) = 2x(m)$. Thus $\alpha(x(m)) = 1$. In the same way we have
for any integer $m$, using the relation $y(n) = 2n + \alpha(n)$ for $n = y(m)$, that
$y(y(m)) = 2y(m) + \alpha(y(m)$, while $y(y(m)) = 2y(m)$. Thus $\alpha(y(m)) = 0$. Since
$X$ and $Y$ form a partition of the integers this gives
$$
n \in X \Leftrightarrow \alpha(n) = 1 \ \ \mbox{\rm and} \ \
n \in Y \Leftrightarrow \alpha(n) = 0.
$$

\medskip

Now we prove that $\alpha(n) = t_2(n)$, i.e., that the sequence $(\alpha(n))_{n \geq 0}$ is the 
Thue-Morse sequence beginning with $0$. It suffices to prove that, for all $m \geq 0$, we have 
$\alpha(2m) = \alpha(m)$ and $\alpha(2m+1) = 1 - \alpha(m)$.

\medskip

If $m$ belongs to $X$, then there exists a $k$ such that $m = x(k)$. We have just seen that 
$\alpha(m) = 1$. We have $x(2m) = 4m + 1 - \alpha(2m)$. But 
$$
x(2m) = x(2x(k)) = x(xx(k)) = xx(x(k)) = 2xx(k) = 4x(k) = 4m.
$$ 
Hence $\alpha(2m) = 1 = \alpha(m)$. Now, since we thus have that $2m$ belongs to $X$, we must have
$2m+1$ belongs to $Y$, hence $\alpha(2m+1) = 0$.

\medskip

If $m$ belongs to $Y$, then there exists a $k$ such that $m = y(k)$. Thus $\alpha(m) = 0$.
We have $y(2m) = 4m + \alpha(2m)$. But
$$
y(2m) = y(2y(k)) = y(yy(k)) = yy(y(k)) = 2yy(k) = 4y(k) = 4m.
$$
Hence $\alpha(2m) = 0$. Now, since we thus have that $2m$ belongs to $Y$, we must have
$2m+1$ belongs to $X$, hence $\alpha(2m+1) = 1$.

\medskip

Finally we thus have that $(\alpha(n))_{n \geq 0} = (t_2(n))_{n \geq 0}$, and then
${\mathbf x} = {\mathbf a}$ and ${\mathbf y} = {\mathbf b}$. \endpf

\bigskip

The next two theorems can be seen as variations on Theorem~\ref{equat}.

\begin{theorem}\label{equat-autre}
Let ${\mathbf x} = (x(n))_{n \geq 0}$ and ${\mathbf y} = (y(n))_{n \geq 0}$ be increasing integer
sequences such that $\{x(n), \ n \geq 0\} \cup \{y(n), \ n \geq 0\} = {\mathbb N}$ satisfying 
$x(0) = 1$, $y(0)=0$, and
$$
\forall n \geq 0, \ x(x(n)) = y(x(n)) - 1 \ \ \mbox{\rm and} \ \ y(y(n)) = x(y(n)) - 1.
$$
Then ${\mathbf x}$ and ${\mathbf y}$ are respectively equal to ${\mathbf a}$ and ${\mathbf b}$ the
sequences of odious and of evil numbers.
\end{theorem}

\proof Let $X = \{x(n), \ n \in {\mathbb N}\}$ and $Y = \{y(n), \ n \in {\mathbb N}\}$. The condition
on $x(x(n))$ and $y(y(n))$ can be written as follows
$$
\mbox{\rm if $m \in X$, then} \ x(m) = y(m) - 1; \
\mbox{\rm if $m \in Y$, then} \ y(m) = x(m) - 1.
$$
This implies in particular that $X \cap Y = \emptyset$, thus $X$ and $Y$ form a partition of the integers.
Now let $1_X$ be the characteristic function of $X$ (i.e., $1_X(n) = 1$ if and only if $n$ belongs to $X$).
Thus $1 - 1_X$ is the characteristic function of $Y$. We will prove by induction on $n$ that, 
for all $n \geq 0$,
$$
x(n) = 2n + 1 - 1_X(n), \ y(n) = 2n + 1_X(n).
$$
The property is true for $n = 0$ since $x(0)=1$ and $y(0)=0$.
If it is true up to $n$, we first have $\{x(n), y(n)\} = \{2n, 2n+1\}$. 
\begin{itemize}

\item If $n+1$ belongs to $X$, then on one hand $x(n+1) = y(n+1)-1$. Since $x(n+1) > \max\{2n, 2n+1\}$,
this gives $x(n+1) \geq 2n+2$ and $y(n+1) = x(n+1)+1 \geq 2n+3$. But $X$ and $Y$ form a partition of
the integers, thus $x(n+1)$ must be equal to $2n+2$ (otherwise $2n+2$ is missed both by $x$ and by $y$),
and $y(n+1) = x(n+1)+1 = 2n+3$. This gives $x(n+1) = 2(n+1) + 1 - 1_X(n+1)$ and $y(n+1) = 2(n+1) + 1_X(n+1)$.

\item If $n+1$ belongs to $Y$, then one one hand $y(n+1) = x(n+1)-1$. Since $y(n+1) > \max\{2n, 2n+1\}$,
this gives $y(n+1) \geq 2n+2$ and $x(n+1) = y(n+1)+1 \geq 2n+3$. But $X$ and $Y$ form a partition of
the integers, thus $y(n+1)$ must be equal to $2n+2$ (otherwise $2n+2$ is missed both by $y$ and by $x$),
and $x(n+1) = y(n+1)+1 = 2n+3$. This gives $y(n+1) = 2(n+1) + 1_X(n+1)$ and $x(n+1) = 2(n+1) + 1 - 1_X(n+1)$.
\end{itemize}

\noindent
We then note that $x(n) = 2n + 1 - 1_X(n)$ for all $n$, implies that $x(x(n)) = 2x(n) + 1 - 1_X(x(n)) = 2n$ 
for all $n$. Similarly $y(n) = 2n + 1_X(n)$ for all $n$ implies that $y(y(n)) = 2y(n)$ for all $n$.
But we have seen that according to $m$ being in $X$ or $Y$, we have $y(m)-x(m) = \pm 1$, i.e.,
$|x(m)-y(m)| = 1$. We can then conclude using Theorem~\ref{equat}. \endpf

\begin{theorem}\label{equat-autre-autre}
Let ${\mathbf x} = (x_n)_{n \geq 0}$ and ${\mathbf y} = (y_n)_{n \geq 0}$ be 
two sequences of integers defined by $x(0) = 1$, $y(0)=0$, and for each $n \geq 1$,
$x(n)$ and $y(n)$ are the smallest integers that did not occur before (i.e., that
do not belong to $\{x(k), \ k \leq n-1\} \cup \{y(k), \ k \leq n-1\}$), with the conditions
that for all $n \geq 0$

\begin{itemize}

\item  $x(x(n))$ and $y(y(n))$ are even,

\item  $x(y(n))$ and $y(x(n))$ are odd.

\end{itemize}

Then ${\mathbf x} = {\mathbf a}$ the sequence of odious numbers, and 
${\mathbf y} = {\mathbf b}$ the sequence of evil numbers.
\end{theorem}

\proof The hypothesis ``the smallest numbers that did not occur before'' implies that $x$ and $y$ do not 
miss any integer. In other words, defining $X = \{x(n), \ n \geq 0\}$ and $Y = \{y(n), \ n \geq 0\}$, we
have $X \cup Y = {\mathbb N}$. On the other hand the intersection of $X$ and $Y$ is empty: if $n$ belongs
both to $X$ and $Y$, then there exist $k, \ell$ with $n = x(k) = y(\ell)$. But then $x(n) = x(x(k))$ is even,
while $x(n) = (x(y(\ell))$ should be odd, which is impossible. Thus $X$ and $Y$ form a partition of the 
integers. We will prove as above that, letting $1_X$ denote the characteristic function of $X$, then
$$
x(n) = 2n + 1 - 1_X(n) \ \mbox{\rm and} \ y(n) = 2n + 1_X(n).
$$
The property is true for $n=0$. Suppose that it is true up to $n$, which implies in particular
that $\{x(n), y(n)\} = \{2n, 2n+1\}$.
\begin{itemize}

\item If $n+1$ belongs to $X$, i.e., $n+1=x(k)$ for some $k$, then $x(n+1) = x(x(k))$ must be even, while 
$y(n+1) = y(x(k)$ must be odd. All integer values being taken, this implies that $x(n+1)=2n+2$ and 
$y(n+1) = 2n+3$. This can also be written $x(n+1) = 2(n+1) + 1 - 1_X(n+1)$ and $y(n+1) = 2(n+1) + 1_X(n+1)$.

\item If $n+1$ belongs to $Y$, i.e., $n+1=y(k)$ for some $k$, then $x(n+1) = x(y(k))$ must be odd, while
$y(n+1) = y(y(k))$ must be even. All integer values being taken, this implies that $y(n+1)=2n+2$ and 
$x(n+1) = 2n+3$. This can also be written $y(n+1) = 2(n+1) + 1_X(n+1)$ and $x(n+1) = 2(n+1) + 1 - 1_X(n+1)$.

\end{itemize}

\noindent
Since for all $n$ we clearly have $|x(n)-y(n)| = 1$ we conclude as in Theorem~\ref{equat-autre}. \endpf

\section{Conclusion}

We would like to add that all the functional equations given above for the sequences of odious and of 
evil numbers can be translated in terms of characterizations of the Thue-Morse sequence. Furthermore
analogous results can be proven for the sequences ${\mathbf a}_{d,j}$.


\begin{thebibliography}{99}

\bibitem{ARS} J.-P. Allouche, N. Rampersad, J. Shallit, {\em On integer sequences 
whose first iterates are linear}, Aequationes Math. {\bf 69} (2005) 114--127.

\bibitem{ubiq} J.-P. Allouche, J. Shallit, {\em The ubiquitous
Prouhet-Thue-Morse sequence}, in Sequences and their Applications
(Singapore, 1998), 1--16, Springer Ser. Discrete Math. Theor. Comput. Sci.,
Springer, London, 1999.

\bibitem{BCK} E. Berlekamp, J. Conway, R. Guy, Winning Ways for Your Mathematical Plays,
vol.~1 and 2, Academic Press, 1982.

\bibitem{BCK2} E. Berlekamp, J. Conway, R. Guy, Winning Ways for Your Mathematical Plays,
vol.~3, Second Edition, A. K. Peters, 2003.

\bibitem {CSV} B. Cloitre, N. J. A. Sloane, M. J. Vandermast, {\em Numerical analogues of 
Aronson's sequence}, J. Integer Seq. {\bf 6} (2003) Art. 03.2.2.

\bibitem{Emmer} M. Emmer, {\em Numeri Malefici (Evil Numbers): Homage to Fabio Mauri}, in 
Imagine Math 2: Between Culture and Mathematics, M. Emmer ed., Springer, 2013, pp. 139--150.

\bibitem{guy} R. K. Guy, {\em Impartial games}, in Combinatorial Games, MSRI Publications, 
Vol. 29, 1995. Available at {\tt http://library.msri.org/books/Book29/files/imp.pdf}

\bibitem{LY} V. Laohakosol, B. Yuttanan, {\em Iterates of increasing sequences of positive 
integers}, Aequationes Math. {\bf 87} (2014) 89--103.

\bibitem{Nietzsche} F. W. Nietzche, Beyond Good and Evil: Prelude to a Philosophy of the Future,
Edited by Rolf-Peter Horstmann, Edited and translated by Judith Norman, Cambridge Texts in the
History of Philosophy, Cambridge, 2002.

\bibitem{oeis} On-Line Encyclopedia of Integer Sequences, available electronically
at the URL {\tt http://oeis.org}

\bibitem{Sarkaria} Karambir S. Sarkaria, {\em Roots of translations}, Aequationes Math. 
{\bf 75} (2008) 304--307. 

\bibitem{Targonski} G. Targo\'nski, Topics in iteration theory, Studia Mathematica: Skript, 6, 
Vandenhoeck \& Ruprecht, G\"ottingen, 1981. 

\bibitem{Wolfram} Wolfram Mathworld. Available at {\tt http://mathworld.wolfram.com/EvilNumber.html}

\end{thebibliography}
\end{document}